# NONPARAMETRIC EMPIRICAL BAYES AND COMPOUND DECISION APPROACHES TO ESTIMATION OF A HIGH-DIMENSIONAL VECTOR OF NORMAL MEANS

By Lawrence D. Brown[1] and Eitan Greenshtein

*University of Pennsylvania and Duke University*

We consider the classical problem of estimating a vector $\boldsymbol{\mu} = (\mu_1, \ldots, \mu_n)$ based on independent observations $Y_i \sim N(\mu_i, 1)$, $i = 1, \ldots, n$.

Suppose $\mu_i$, $i = 1, \ldots, n$ are independent realizations from a completely unknown $G$. We suggest an easily computed estimator $\hat{\boldsymbol{\mu}}$, such that the ratio of its risk $E(\hat{\boldsymbol{\mu}} - \boldsymbol{\mu})^2$ with that of the Bayes procedure approaches 1. A related compound decision result is also obtained.

Our asymptotics is of a triangular array; that is, we allow the distribution $G$ to depend on $n$. Thus, our theoretical asymptotic results are also meaningful in situations where the vector $\boldsymbol{\mu}$ is sparse and the proportion of zero coordinates approaches 1.

We demonstrate the performance of our estimator in simulations, emphasizing sparse setups. In "moderately-sparse" situations, our procedure performs very well compared to known procedures tailored for sparse setups. It also adapts well to nonsparse situations.

**1. Introduction.** Let $Y = (Y_1, \ldots, Y_n)$ be a random normal vector where $Y_i \sim N(\mu_i, 1)$, $i = 1, \ldots, n$ are independent. Consider the classical problem of estimating the mean vector $\boldsymbol{\mu} = (\mu_1, \ldots, \mu_n)$ by a (nonrandomized) estimator $\Delta = \Delta(Y)$ under the squared-error loss $L_n(\boldsymbol{\mu}, \Delta) = \sum_i (\Delta_i - \mu_i)^2$. The corresponding risk function is the expected squared error

$$R(\boldsymbol{\mu}, \Delta) = E_{\boldsymbol{\mu}}(L_n(\boldsymbol{\mu}, \Delta(Y))).$$

*Compound decision theory.* A natural class of decision functions is the class of simple symmetric estimators that was suggested by Robbins (1956). This is the class of all estimators $\Delta^*$ of the form

$$\Delta^*(Y) = (\delta(Y_1), \ldots, \delta(Y_n))$$

---

Received December 2007; revised May 2008.
[1]Supported in part by NSF Grant DMS-07-07033.

for some function, $\delta$. For such an estimator, we will occasionally write $\Delta^*(Y) = \Delta^*(Y|\delta)$ in order to show the dependence on $\delta$.

Given $\boldsymbol{\mu} = (\mu_1, \ldots, \mu_n)$, let

$$\delta^{*\boldsymbol{\mu}} = \arg\min_\delta R(\boldsymbol{\mu}, \Delta^*(\cdot|\delta))$$

and, for notational convenience, let $\Delta^{*\boldsymbol{\mu}} = \Delta(\cdot|\delta^{*\boldsymbol{\mu}})$.

Consider an oracle that knows the value of the vector $\boldsymbol{\mu}$ but must use a simple-symmetric estimator. Such an oracle would use the estimator $\Delta^{*\boldsymbol{\mu}}$. The goal of compound decision theory is to achieve *nearly* the risk obtained by such an oracle, but by using a "legitimate" estimator, one that may involve the entire vector of observations Y but does not involve knowledge of the parameter vector $\boldsymbol{\mu}$. In establishing specific results, it is important to be suitably precise about the (asymptotic) sense in which this near-ness is measured. This will be discussed later, after introducing the companion concept of empirical Bayes, for background on both compound decision and empirical Bayes [see Robbins (1951, 1956, 1964), Samuel (1965), Copas (1969) and Zhang (2003), among many other papers].

*Empirical Bayes.* Let $G$ be a prior distribution on $\mathcal{R}$. Let $\mathbf{M} = \{M_i, i = 1, \ldots, n\}$ be an unobserved random sample from this distribution. Conditional on the $\{M_i\}$ observe $Y_i \sim N(M_i, 1), i = 1, \ldots, n$, independent. Here, the target procedure is the Bayes procedure, to be denoted $\Delta^G$. The goal is to find a procedure $\Delta$ whose expected risk under $G$ is suitably near that of $\Delta^G$ as $n \to \infty$, when $G$ is unknown. The notation here is intentionally similar to that used previously for the compound decision problem, but note that the superscript is now a distribution $G$, whereas, in the compound decision situation, the superscript is a vector $\boldsymbol{\mu}$ or, equivalently, the set of coordinates of $\boldsymbol{\mu}$.

*Relation of compound and empirical Bayes risks.* The expected average risk under $G$ of a procedure $\Delta$ will be denoted by $B(G, \Delta)$. Note that

$$(1) \qquad B(G, \Delta) = E_G\left(\frac{1}{n} R(\mathbf{M}, \Delta)\right),$$

where we treat $\mathbf{M}$ as a random vector whose coordinates are a sample of size $n$ from $G$, as described above. (For convenience, the dependence on $n$ is suppressed in the notation.)

Here are some simple consequences of this relation. Let $\{\Delta_n\}$ denote a sequence of estimators in a sequence of problems with increasing dimension $n$. Suppose, for example, that $\{\Delta_n\}$ has the basic asymptotic compound Bayes property that, for every $\boldsymbol{\mu}^n = (\mu_1^n, \ldots, \mu_n^n)$, $n = 1, 2, \ldots$,

$$(2) \qquad \frac{1}{n} R(\boldsymbol{\mu}^n, \Delta_n) - \frac{1}{n} R(\boldsymbol{\mu}^n, \Delta^{*\boldsymbol{\mu}^n}) \overset{n \to \infty}{\to} 0.$$

Also, assume that $\frac{1}{n}R(\boldsymbol{\mu}, \Delta_n)$ is uniformly bounded, as will typically be the case under suitable assumptions [as in (35)]. It then follows from (1) that $\{\Delta_n\}$ is asymptotically empirical Bayes in the basic sense that

$$\begin{aligned}
B(G, \Delta^G) &= E_G\left(\frac{1}{n}R(\mathbf{M}, \Delta^G)\right) \geq E_G \frac{1}{n}(R(\mathbf{M}, \Delta^{*\mathbf{M}})) \\
&= E_G\left(\frac{1}{n}R(\mathbf{M}, \Delta_n)\right) + o(1) = B(G, \Delta_n) + o(1).
\end{aligned} \tag{3}$$

Hence, under very mild conditions, asymptotic compound optimal in the sense of (2), implies asymptotic empirical Bayes in the sense of (3).

In Section 2, we will propose a particular, easily implemented form for $\Delta_n$. In Section 3, we establish some more precise compound and empirical Bayes properties for this estimator. Although these properties are more demanding than (2) and (3), the relation (1) remains an important part of the arguments that establish them.

*Relation of compound optimal and empirical Bayes procedures.* The relation (1) describes a close connection between the compound Bayes and empirical Bayes criteria. It is also true that the optimal procedures are closely connected. The Bayes procedure $\Delta^G(Y) = (\delta^G(Y_1), \ldots, \delta^G(Y_n))$, for a specified prior $G$, is of course given by Bayes formula

$$\delta^G(Y_j) = E(M_j|Y) = E(M_j|Y_j) = \frac{\int u\phi(u - Y_j)G(du)}{\int \phi(u - Y_j)G(du)}. \tag{4}$$

Among its other features, the Bayes formula (4) reveals that the Bayes procedure is a simple-symmetric estimator.

A simple derivation also yields the basic formula for $\Delta^{*\boldsymbol{\mu}}$, through the corresponding univariate function $\delta^{*\boldsymbol{\mu}}$

$$\delta^{*\boldsymbol{\mu}}(u) = \frac{\sum_i \mu_i \phi(\mu_i - u)}{\sum_i \phi(\mu_i - u)}. \tag{5}$$

Given $\boldsymbol{\mu} = (\mu_1, \ldots, \mu_n)$, let $F_n^{\boldsymbol{\mu}}$ denote the corresponding empirical CDF. Then the formula for $\delta^{*\boldsymbol{\mu}}$ can be rewritten as

$$\delta^{*\boldsymbol{\mu}}(u) = \delta^{F_n^{\boldsymbol{\mu}}}(u). \tag{6}$$

This formula provides a direct connection between the optimal estimators for the two settings. This will be exploited in the construction, in Section 2, of an asymptotically optimal estimator.

*Sparse estimation problems.* The following discussion is intended to help motivate the asymptotic properties to be established in Section 3. It will also help motivate the choice of settings serving as the basis for the numerical results reported in Section 4.

Many recent statistical results have focused on the importance of treating situations involving "sparse" models [see, e.g., Donoho and Johnstone (1994), Johnstone and Silverman (2004) and Efron (2003)]. Many such problems involve issues of testing hypotheses, but for others estimation is of secondary or even primary interest. The basic asymptotic empirical Bayes property in (3) involves asymptotic properties for a fixed prior $G$. Such a formulation is not sufficiently flexible to provide useful results in "sparse" settings. In Section 3, we investigate an asymptotic formulation that is appropriate for many sparse problems, as well as for the more conventional settings involving asymptotics for fixed (but unknown) $G$.

"Sparsity" is not a precise statistical condition. However, the essence of many "sparse" settings is captured by considering situations in which most of the unknown coordinates $\mu_i$ take the value 0, and the remaining few take other value(s).

To be precise, in the following discussion of the compound Bayes setting, consider a situation in which the possible values for the coordinates $\mu_i \equiv \mu_i^n$ of $\boldsymbol{\mu}^n$ are either $\mu_i = 0$ or $\mu_i = \mu_0 \neq 0$, $i = 1, \ldots, n$. Here, we consider a sequence of problems with increasing dimension $n$. For a given $n$, let $p = p(n)$ denote the proportion of nonzero values. The situation is *sparse* if $p(n) \to 0$ as $n \to 0$. (For simplicity, assume that there is only one possible nonzero value, $\mu_0$, and that this value does not change with $n$. Of course many other situations are possible that should still be classed as sparse models.) Then,

$$(7) \qquad \frac{1}{n} R(\boldsymbol{\mu}^n, \Delta^{*\boldsymbol{\mu}^n}) = O(p(n)).$$

Note that

$$(8) \qquad \frac{1}{n} R(\boldsymbol{\mu}_n, \Delta^{*\boldsymbol{\mu}^n}) \stackrel{p(n) \to 0}{\to} 0.$$

Hence, useful asymptotic results for sparse models must accommodate this fact.

The asymptotic statements in Section 3 are naturally scaled to accommodate sparsity in this way because they examine the relative risk ratio, rather than the ordinary difference between average risks, as in the basic statement (2). Thus, for the given sequence, $\{\hat{\Delta}_n\}$ of procedures defined in Section 2, these results examine the limiting value of

$$(9) \qquad \frac{R(\boldsymbol{\mu}^n, \hat{\Delta}_n) - R(\boldsymbol{\mu}^n, \Delta^{*\boldsymbol{\mu}^n})}{R(\boldsymbol{\mu}^n, \Delta^{*\boldsymbol{\mu}^n})}$$

and establish quite general conditions under which this ratio converges to 0. (The results of Section 3 include the preceding two point model as a very special case.)

Here is the empirical Bayes setting which corresponds to the special sparse compound Bayes model described in the previous paragraphs. Consider an empirical Bayes model, in which it is assumed that $G = G_n$, where

$$G_n(\{\mu_0\}) = \pi_i(n) = 1 - G_n(\{0\}). \tag{10}$$

Note that, as in (7),

$$B(G_n, \Delta^{G_n}) = O(\pi(n)).$$

Similar to (9), the asymptotic results appropriate for sparse models will be phrased in terms of the limiting value of the ratio difference

$$\frac{B(G_n, \hat{\Delta}_n) - B(G_n, \Delta^{G_n})}{B(G_n, \Delta^{G_n})}. \tag{11}$$

The preceding discussion suggests that the degree of "sparsity" of a sequence of compound or empirical Bayes models could be measured by the asymptotic behavior of $R(\boldsymbol{\mu}^n, \Delta^{*\boldsymbol{\mu}^n})$ or $B(G_n, \Delta^{G_n})$, respectively. For example, sequences of models for which

$$\liminf \frac{1}{n} R(\boldsymbol{\mu}^n, \Delta^{*\boldsymbol{\mu}^n}) > 0 \tag{12}$$

[or $\liminf B(G_n, \Delta^{G_n}) > 0$] could be considered *nonsparse*. At the other extreme are sequences for which

$$R(\boldsymbol{\mu}^n, \Delta^{*\boldsymbol{\mu}^n}) = O(1); \tag{13}$$

those could be called extremely sparse. Sequences between those extremes can be termed moderately-sparse. A typical example could be a sequence of problems for which

$$R(\boldsymbol{\mu}^n, \Delta^{*\boldsymbol{\mu}^n}) = O(n^\alpha), \qquad 0 < \alpha < 1. \tag{14}$$

Note that, in this description of sparseness, the zero value does not play a special role. It is the "complexity" of the sequence or the "difficulty to estimate it" that defines its sparseness.

In Section 2, we construct an estimator that is approximately compound optimal and empirical Bayes. The construction formula is simple and easily implemented. This estimator performs very well for nonsparse and moderately sparse settings, such as those in (14). It can also be satisfactorily used for extremely sparse settings, but it is implicit in the theory in Section 3 and explicit in the simulations in Section 4 that its performance is not quite optimal in some extremely sparse settings.

*Permutation invariant procedures.* A natural class of procedures, which is larger than the class of simple symmetric ones, is the class of permutation invariant procedures. This is the class of all procedures $\Delta$ that satisfy

$$\Delta(Y_1,\ldots,Y_n) = (\hat{\mu}_1,\ldots,\hat{\mu}_n) \quad \Leftrightarrow \quad \Delta(Y_{\pi(1)},\ldots,Y_{\pi(n)}) = (\hat{\mu}_{\pi(1)},\ldots,\hat{\mu}_{\pi(n)})$$

for every permutation $\pi$.

In a recent paper by Greenshtein and Ritov (2008), a "strong equivalence" between the class of permutation invariant procedures and the class of simple symmetric procedures is shown. This equivalence implies that some of the optimality results we obtain, comparing the performance of our procedure with that of the optimal simple symmetric procedure for a given $\boldsymbol{\mu}$, are valid also with respect to the comparison with the (superior) optimal permutation invariant procedure.

**2. Bayes, empirical Bayes and compound decision.** Let $Y \sim N(M,1)$ where $M \sim G$, $G \in \mathcal{G}$. We want to emulate the Bayes procedure $\delta^G \equiv \delta_1^G$, based on a sample $Y_1,\ldots,Y_n$, $Y_i \sim N(M_i,1)$, $i=1,\ldots,n$, where $M_i \sim G$ and the $Y_i$ are independent conditional on $M_1,\ldots,M_n$, $i=1,\ldots,n$. In general, $G$ may depend on $n$, but, in order to simplify the notation and presentation, we consider a fixed $G$ throughout this section. The generalization for a triangular array is easily accomplished.

Consider our problem for a general variance $\sigma^2$; that is, suppose $Y_i \sim N(M_i, \sigma^2)$, $M_i \sim G$, $i = 1,\ldots,n$. Let $g_{G,\sigma^2}^*$ be the mixture density

$$(15) \qquad g_{G,\sigma^2}^*(y) = \int \frac{1}{\sigma} \phi\left(\frac{y-\mu}{\sigma}\right) dG(\mu).$$

Then, from Brown (1971), (1.2.2), we have that the Bayes procedure denoted $\delta_{\sigma^2}^G$, satisfies

$$(16) \qquad \delta_{\sigma^2}^G(y) = y + \sigma^2 \frac{g_{G,\sigma^2}^{*\prime}(y)}{g_{G,\sigma^2}^*(y)}.$$

Here, $g_{G,\sigma^2}^{*\prime}(y)$ is the derivative of $g_{G,\sigma^2}^*(y)$.

The estimator that we suggest for $\delta_1^G$ is of the form

$$(17) \qquad \hat{\delta} = y + \frac{\hat{g}_h^{*\prime}(y)}{\hat{g}_h^*(y)},$$

where $\hat{g}_h^{*\prime}(y)$ and $\hat{g}_h^*(y)$ are appropriate kernel estimators for the density $g_{G,1}^*(y)$ and its derivative $g_{G,1}^{*\prime}(y)$. The subscript $h$ is the bandwidth for the estimator. We will use a normal kernel. This choice is convenient from several perspectives, but does not seem to be essential. See Remark 2 later in this section. An alternative to kernel density estimators could be a direct estimation of $G$. An approach involving MLE estimation of $G$ was recently suggested by Wenhua and Zhang (2007). Its performance in simulations is excellent and it has appealing theoretical properties. However, it is computationally intensive.

Let $h > 0$ be a bandwidth constant. Typically, $h$ will depend on $n$, and $\lim_{n\to\infty} h = 0$. Then, define the kernel estimator

$$\hat{g}_h^*(y) = \frac{1}{nh} \sum \phi\left(\frac{y - Y_i}{h}\right). \tag{18}$$

Its derivative has the form

$$\hat{g}_h^{*\prime}(y) = \frac{1}{nh} \sum \frac{Y_i - y}{h^2} \times \phi\left(\frac{y - Y_i}{h}\right). \tag{19}$$

Let

$$v = 1 + h^2.$$

The following simple lemma establishes that $\hat{g}_h^*$ and $\hat{g}_h^{*\prime}$ are unbiased estimates of $g_{G,v}^*$ and $g_{G,v}^{*\prime}$. It also further interprets their form.

Let $G_n^Y$ denote the empirical distribution determined by $Y_1, \ldots, Y_n$.

LEMMA 1. *Let $h > 0$ and $v = 1 + h^2$, and suppose $Y_i \sim N(M_i, 1)$, where $M_i \sim G$ are independent. Then,*

$$\hat{g}_h^*(y) = g_{G_n^Y, v-1}^*(y), \qquad \hat{g}_h^{*\prime}(y) = g_{G_n^Y, v-1}^{*\prime}(y), \tag{20}$$

$$E g_{G_n^Y, v-1}^*(y) = g_{G,v}^*(y), \qquad E g_{G_n^Y, v-1}^{*\prime}(y) = g_{G,v}^{*\prime}(y). \tag{21}$$

PROOF. We write

$$\hat{g}_h^*(y) = \int \frac{1}{h} \phi\left(\frac{y-t}{h}\right) dG_n^Y(dt) = g_{G_n^Y, h^2}^*(y),$$

since $h^{-1}\phi(x/h)$ is the normal density with variance $h^2$. Let $\Phi_{\sigma^2}$ denote the normal distribution with variance $\sigma^2$. Under the conditions of the lemma, $E(G_n^Y) = G * \Phi_1$.

Hence, $E(\hat{g}_h^*(y)) = g_{G*\Phi_1, h^2}^*(y) = g_{G, 1+h^2}^*(y)$, since $(G * \Phi_1) * \Phi_{h^2} = G * \Phi_{1+h^2}$.

The arguments for the derivatives follow by differentiation or by an independent argument analogous to the above. This completes the proof. □

Hence, the basic formula (17) may be rewritten as

$$\hat{\delta}_{1+h^2}(y) = \hat{\delta}_v(y) = y + \frac{g_{G_n^Y, h^2}^{*\prime}(y)}{g_{G_n^Y, h^2}^*(y)}. \tag{22}$$

As a final step in the motivation of our estimator, note that

(23) $$\delta^G_{1+h^2}(y) = \delta^G_v(y) \stackrel{h\to 0}{\Rightarrow} \delta^G_1(y).$$

By Lemma 1 and (23), we expect that, for large $n$ and $v = 1 + h^2 \approx 1$, we have

(24) $$\frac{g^{*\prime}_{G^Y_n, v-1}(y)}{g^*_{G^Y_n, v-1}(y)} \approx \frac{g^{*\prime}_{G,v}(y)}{g^*_{G,v}(y)} \approx \frac{g^{*\prime}_{G,1}(y)}{g^*_{G,1}(y)}.$$

Similarly, we have

$$\delta^G_1(y) = y + (\delta^G_1(y) - y) \approx y + \left(\left[y + \frac{g^{*\prime}_{G,v}(y)}{g^*_{G,v}(y)}\right] - y\right)$$

(25) $$\approx y + \frac{1}{v-1}\left(\left[y + (v-1)\frac{g^{*\prime}_{G^Y_n, v-1}(y)}{g^*_{G^Y_n, v-1}(y)}\right] - y\right)$$

$$= y + \frac{1}{v-1}(\delta^{G^Y_n}_{v-1}(y) - y) = \hat{\delta}_v(y).$$

Here, $\delta^{G^Y_n}_{v-1}$ is as defined above (16).

REMARK 1. All the equations obtained so far for the empirical Bayes setup have a parallel derivation and presentation in the compound decision setup for a given $\boldsymbol{\mu} = (\mu_1, \ldots, \mu_n)$, where $F^{\boldsymbol{\mu}}_n$, the empirical distribution of $\boldsymbol{\mu} = (\mu_1, \ldots, \mu_n)$, plays the role of $G$, as in (6). For example, (15) has the form $\frac{1}{n}\sum_i \frac{1}{\sigma}\phi(\frac{y-\mu_i}{\sigma})$, and the analog of $\delta^G_{\sigma^2}$ is denoted $\delta^{*\boldsymbol{\mu}}_{\sigma^2}$, etc.

EXAMPLE 1. It is of some interest to examine how the preceding formulas compare in the standard case where the true prior is Gaussian, say $G \sim N(0, \gamma^2)$. In that case, $G^Y_n \Rightarrow N(0, 1+\gamma^2)$ in distribution. The actual Bayes procedure is

$$\delta^G_1(Y) = \left(1 - \frac{1}{1+\gamma^2}\right)Y,$$

while, by (25) for a fixed $v$, $\hat{\delta}_v(Y)$ converges as $n \to \infty$ to

$$\left(1 - \frac{1}{v+\gamma^2}\right)Y.$$

This may be seen when substituting $\delta^{N(0,1+\gamma^2)}_{v-1}(y) = (1 - \frac{v-1}{v+\gamma^2})y$, for $\delta^{G^Y_n}_{v-1}(y)$ in (25).

Thus, when letting $v \equiv v_n$ approach 1 (equivalently when letting the bandwidth $h = \sqrt{v-1}$ approach 0) as $n$ approaches infinity, we may see that $\hat{\delta}_v(y)$ approaches $\delta^G_1(y)$.

REMARK 2 (*On the choice of a kernel*). One could choose other kernels and obtain corresponding different estimators. See, for example, the papers of Zhang (1997, 2005). In those papers, Zhang introduces an estimator for $\delta^G$, using Fourier methods and corresponding kernels, to estimate the above $g_{G,1}^{*\prime}$ and $g_{G,1}^*$. Zhang's papers are very relevant, and there are similarities between our approach and his earlier development.

We now point to some advantages of our kernel. One advantage is the interpretation of $\hat{\delta}_v$ as an approximation for $\delta_v^G$. Here, $\delta_v^G(u)$ is the Bayes decision function for the setup where $U \sim N(M, v)$ and $M \sim G$ (see, e.g., Example 1, with the interpretation of the obtained rules, in terms of the approximation $G_n^Y$ of $G$). This interpretation is very helpful in the proof of Theorem 1. We are not sure to what extent a normal kernel is essential to obtain the good performance of our estimator, but it certainly simplifies various arguments. In addition, kernels with heavy tails would typically introduce a significant bias when estimating $g_{G,1}^*$ and its derivative in the tail.

**3. Optimality in compound decision under sparsity.** In this section, we study asymptotics which are appropriate for both nonsparse and sparse compound decision problems. The traditional asymptotics for empirical Bayes and compound decision, consider the difference in average risks between the target (or optimal) procedure and a suggested estimator. In the sparse setting, both of these quantities approach zero. So the traditional asymptotic criteria are not informative, and a more delicate study is needed.

Our main result, Theorem 1, covers the compound decision framework. It has an analogous empirical Bayes formulation which is obtained as a corollary.

The formal setup is of a triangular array, where, at stage $n$, the parameter space, denoted $\Theta^n$, is of dimension $n$. We use the notation $\boldsymbol{\mu}^n = (\mu_1^n, \ldots, \mu_n^n) \in \Theta^n$.

For every $\varepsilon > 0$, we assume

(26) $\qquad |\mu_j^n| < C_n = o(n^\varepsilon), \qquad n = 1, 2, \ldots, \infty, \ j = 1, \ldots, n.$

Such configurations include the interesting cases where $\mu_j^n = O(\sqrt{\log(n)})$. Those are interesting configurations in which the statistical task of discriminating between signal and noise is neither too easy nor too hard.

As before, we observe a vector $(Y_1^n, \ldots, Y_n^n)$, where $Y_j^n$, are independently distributed $N(\mu_j^n, 1)$. Consider the loss for estimating $\boldsymbol{\mu}^n$ by $\hat{\boldsymbol{\mu}}^n$,

(27) $$L(\boldsymbol{\mu}^n, \hat{\boldsymbol{\mu}}^n) = \sum_{j=1}^n (\hat{\mu}_j^n - \mu_j^n)^2,$$

here $\hat{\boldsymbol{\mu}}^n = (\hat{\mu}_1^n, \ldots, \hat{\mu}_n^n)$.

In this section, we will introduce the following slight modification for $\hat{\delta}_v(u)$, and will consider a truncated estimator which at stage $n$ is of the form

$$\hat{\delta}_v^t(u) = \text{sign}(\hat{\delta}_v(u)) \times \min(C_n, |\hat{\delta}_v(u)|). \tag{28}$$

Note that we chose to truncate $\hat{\delta}_v$ so that $|\hat{\delta}_v| < C_n$. An alternate truncation can be used that may be more desirable in practice. This involves truncation of $\hat{\delta}_v(y) - y$, rather than $\hat{\delta}_v$. In this case, the truncation level can be chosen independent of $C_n$. We write $\hat{\delta}_v(y) = y + (\hat{\delta}_v(y) - y) \equiv y + R$. Let $\tilde{R} = \text{sign}(R)\min(|R|, \sqrt{3\log(n)})$. The alternate truncated estimator is

$$y + \tilde{R}. \tag{29}$$

This estimator also satisfies the conclusion of Theorem 1 and our other results. Minor modifications of the proofs are needed. Let $\hat{\Delta}_v^t(Y) = \Delta^*(Y|\hat{\delta}_v^t)$ denote the simple symmetric estimator of $\boldsymbol{\mu}^n$. Recall $v = 1 + h^2$. We now state our main result.

THEOREM 1. *Consider a triangular array with $\Theta^n$, as above, and sequences $\boldsymbol{\mu}^n \in \Theta^n$ as in (26). Let $v \equiv v_n \to 1$, $v > 1$, be any sequence satisfying:*

(i) $\frac{1}{v-1} = o(n^{\varepsilon'})$ *for every* $\varepsilon' > 0$.
(ii) $\log(n) = o(\frac{1}{v-1})$.

*Assume that, for some $\varepsilon > 0$ and $n_0$,*

$$R(\boldsymbol{\mu}^n, \Delta^{\boldsymbol{\mu}^n}) > n^\varepsilon \qquad \forall n > n_0. \tag{30}$$

*Then,*

$$\limsup \frac{R(\boldsymbol{\mu}^n, \hat{\Delta}_v^t)}{R(\boldsymbol{\mu}^n, \Delta^{*\boldsymbol{\mu}^n})} = 1. \tag{31}$$

REMARK 3. Theorem 1 states that, in situations which are not too advantageous for the oracle so that its risk is of an order larger than $n^\varepsilon$ for some $\varepsilon > 0$, we may asymptotically do as well as that oracle by letting $v$ approach 1 in the right way. Doing as well as the oracle means that the ratio of the risks approaches 1. Note that some condition resembling (30) is needed; if, for example, $\boldsymbol{\mu}^n = (0, \ldots, 0)$, $n = 1, 2, \ldots$, then the corresponding risk of the oracle is identically 0, and we can obviously not achieve such a risk by our estimator.

Although the asymptotics in this section are motivated mainly by sparse setups, the result in Theorem 1 is valid for any sequence $\boldsymbol{\mu}^n$ satisfying (26) and (30). Obtaining an estimator that performs well and adapts well to a broad range of "sparseness"/"denseness," is the main achievement in this paper. The simple, easily interpretable form of our estimator is an additional useful feature.

REMARK 4. There is an alternate form for the conclusion (31) that avoids the necessity for an explicit assumption like (30). A minor additional argument shows that in the statement of the theorem one can omit (30) and replace the conclusion (31) by the conclusion

$$\limsup \frac{R(\boldsymbol{\mu}^n, \hat{\Delta}_v^t)}{R(\boldsymbol{\mu}^n, \Delta^{*\boldsymbol{\mu}^n}) + A_n} \leq 1 \tag{32}$$

for all sequences $\{A_n\}$ such that $A_n > n^{\varepsilon'}$ for some $\varepsilon' > 0$.

REMARK 5 (*On the choice of the bandwidth*). The asymptotic result of the theorem requires that $h_n \to 0$, but at a fairly slow rate. This slow rate is needed in order to obtain the general conclusion in (31), assuming any value of $\varepsilon$ in (30). However, when (30), holds for large values of $\varepsilon$ (e.g., nonsparse case with $\varepsilon = 1$), then smaller values of $h_n$ might be desirable and will have some theoretical advantage. Our theoretical results suggest that $h_n^2$ should converge to zero "just faster" than $1/\log(n)$; we recommend $h_n^2 = 1/\log(n)$ as a "practical default choice." This choice was studied in our simulations and also in Brown (2008) and Greenshtein and Park (2007), where real data sets are explored.

One could improve by selecting different values of bandwidths for different points $y$ in an adaptive manner. Obviously, smaller bandwidth are desirable in the "main body" of the distribution and bigger ones on the tail. Also, one could use different bandwidth when estimating the density and its derivative at a point (typically larger bandwidth for estimating the derivative). Such an approach (and possible improvement) would introduce computational complexity to our simply computed estimator. We do not pursue this approach in the present manuscript.

From now on, we will occasionally drop the superscripts $t$ in $\hat{\delta}_v^t$, and $n$ on $\boldsymbol{\mu}^n$. Recall the notation $\delta_v^{*\boldsymbol{\mu}}$ for the optimal simple symmetric function given $\boldsymbol{\mu}$, when $Y_i \sim N(\mu_i, v)$. Thus, $\delta^{*\boldsymbol{\mu}} \equiv \delta_1^{*\boldsymbol{\mu}}$.

Write

$$\sum (\hat{\delta}_v(Y_i) - \mu_i)^2 = \sum (\hat{\delta}_v(Y_i) - \delta_v^{*\boldsymbol{\mu}}(Y_i) + \delta_v^{*\boldsymbol{\mu}}(Y_i) - \mu_i)^2.$$

Theorem 1 will follow when we prove the following two lemmas and apply Cauchy–Schwarz.

LEMMA 2. *For $v \equiv v_n > 1$, such that $\log(n) = o(\frac{1}{v-1})$, and $\boldsymbol{\mu}^n \in \Theta^n$ as in (26),*

$$\lim \frac{E_{\boldsymbol{\mu}^n} \sum (\delta_1^{*\boldsymbol{\mu}}(Y_i^n) - \mu_i^n)^2}{E_{\boldsymbol{\mu}^n} \sum (\delta_v^{*\boldsymbol{\mu}}(Y_i^n) - \mu_i^n)^2} = 1. \tag{33}$$

PROOF. See Appendix. □

LEMMA 3. *Let $\varepsilon > 0$ (arbitrarily small). Suppose that $v \equiv v_n > 1$, satisfy $\frac{1}{v-1} = o(n^{\varepsilon'})$ for every $\varepsilon' > 0$, and $\boldsymbol{\mu}^n \in \Theta^n$ as in (26). Then,*

$$(34) \qquad E_{\boldsymbol{\mu}^n} \sum_i (\delta_v^{*\boldsymbol{\mu}}(Y_i^n) - \hat{\delta}_v^t(Y_i^n))^2 = o(n^\varepsilon).$$

PROOF. See Appendix. □

A result analogous to Theorem 1 for the empirical Bayes setup is obtained as a corollary. Consider a triangular array where at stage $n$, we observe $Y_i^n \sim N(M_i^n, 1)$, $M_i^n \sim G_n$, $i = 1, \ldots, n$, $M_i^n$ are independent and $Y_i^n$ are independent conditional on $M_i^n$, $i = 1, \ldots, n$; $G_n$ are unknown. Assume that $G_n$ have a support on $(-C_n, C_n)$, where

$$(35) \qquad C_n = o(n^{\varepsilon'})$$

for every $\varepsilon' > 0$. Let $\delta_1^{G_n}$ be the sequence of Bayes procedures. In the following corollary, the expectation is taken with respect to the joint distribution of $(M_1^n, Y_1^n), \ldots, (M_n^n, Y_n^n)$.

COROLLARY 1. *Let $\varepsilon > 0$ (arbitrarily small). For any sequence $v = v_n > 1$, such that:*

(i) $\frac{1}{v-1} = o(n^{\varepsilon'})$ for every $\varepsilon' > 0$,
(ii) $\log(n) = o(\frac{1}{v-1})$,

$$(36) \qquad \limsup \frac{E \sum_i (\hat{\delta}_v^t(Y_i^n) - M_i^n)^2}{E \sum (\delta_1^{G_n}(Y_i^n) - M_i^n)^2 + n^\varepsilon} \le 1.$$

PROOF. The corollary is obtained by conditioning on every possible realization $\mathbf{M}^n = (M_1^n, \ldots, M_n^n)$ and applying Theorem 1 coupled with Remark 4 on each realization treating the conditional setup as a compound decision problem. The proof follows, since, by definition, for every $\mathbf{M}^n$ (treated as a fixed vector), $\sum_i E_{M_i^n}(\delta^{*\mathbf{M}^n}(Y_i) - M_i^n)^2 \le \sum_i E_{M_i^n}(\delta_1^{G_n}(Y_i) - M_i^n)^2$. □

REMARK 6. Assuming the more restrictive condition $C_n = \sqrt{K \log(n)}$, for some $K$, a careful adaptation of our proof will yield the conclusion of Theorem 1 under the weaker assumption that $R(\boldsymbol{\mu}^n, \Delta^{\boldsymbol{\mu}^n})$ is a suitable power of $\log(n)$.

**4. Simulations.** This section will demonstrate the performance of our method in a range of settings. As explained, the value of $v$ should decrease as $n$ increases and should be chosen bigger than 1 but close to 1. We used $v = 1.15$ in simulations with $n = 1000$, $v = 1.1$ when $n = 10{,}000$ and $v = 1.05$ when $n = 100{,}000$. No attempt was made to optimize $v$. Note our default recommendation choice, $v = 1 + (1/\log(n))$ equals 1.144 and 1.108 for $n = 1000$ and $n = 10{,}000$, correspondingly, roughly according to our choice. For $n = 100{,}000$, we chose $v = 1.05$ rather than 1.086 in order to keep a gap of 0.05. However, small changes (say, take $v = 1.15$ rather than $v = 1.1$) did not have much of an effect.

In Table 1 of Johnstone and Silverman (2004) [cited bellow as J–S (2004)], the performances of eighteen estimation methods were compared in various setups and configurations. Those methods include soft and hard universal thresholds and others. The performance was compared in terms of the expected squared risk. In all the configurations, the dimension of the vector $\boldsymbol{\mu}$ is $n = 1000$ [i.e., $Y_i \sim N(\mu_i, 1)$, $i = 1, \ldots, 1000$]. In four configurations, there are $k = 5$ nonzero signals and these nonzero signals all take the value $u_1 = 3$, or all are $u_1 = 4$, $u_1 = 5$ or $u_1 = 7$, respectively. A similar study was done when there are $k = 50$ nonzero signals, and $k = 500$ nonzero signals; the values of the nonzero signals are as before.

In the second line of the following Table 1, we show the performance of the best among the eighteen methods in each case (i.e, the performance of the method with minimal simulated risk for the specific configuration). The first line shows the performance of our $\tilde{\delta}_v$ with $v = 1.15$. The performance and empirical risk of our procedure is based on averaging of the results of 50 simulations in each configuration. One can see that the empirical risk of our procedure is lower than the minimum of all the others in the nonsparse case and in the moderately sparse case. Our procedure adapts particularly well in the nonsparse case. Our method does not do that well in the extremely sparse case; it is worse than the various empirical Bayes procedures suggested in Johnstone and Silverman's paper, but it is within the range of the other methods. All entries, in this table and those to follow, are rounded to the nearest integer.

TABLE 1
*Risk of $\tilde{\delta}_{1.25}$ compared to that of the best procedure in J–S (2004); $n = 1000$ (average of 50 simulations rounded to the nearest integer)*

| $k$ | 5 | | | | 50 | | | | 500 | | | |
|---|---|---|---|---|---|---|---|---|---|---|---|---|
| $u_1$ | 3 | 4 | 5 | 7 | 3 | 4 | 5 | 7 | 3 | 4 | 5 | 7 |
| $\tilde{\delta}_{1.15}$ | 53 | 49 | 42 | 27 | 179 | 136 | 81 | 40 | 484 | 302 | 158 | 48 |
| Minimum | 34 | 32 | 17 | 7 | 201 | 156 | 95 | 52 | 829 | 730 | 609 | 505 |

Note that, when the risk of the oracle is very small, our Theorem 1 does not imply that we are doing well with respect to the oracle in terms of risk ratios.

Here, $\tilde{\delta}_v$ denotes the following minor adaptation of $\hat{\delta}_v$:

$$\tilde{\delta}_v = y + v \frac{\hat{g}^{*\prime}_{G_n,v-1}(y)}{\hat{g}^{*}_{G_n,v-1}(y)}.$$

The difference, relative to $\hat{\delta}_v$, is the multiplication by $v$ of the ratio. As $v \to 1$ the difference between the two estimators is negligible. The procedure $\tilde{\delta}_v$ seems more suitable in approximating, $\delta^{*\mu}_v$ and is as appealing as $\hat{\delta}_v$.

In the following Table 2, we report on the behavior of our procedure based on 50 simulations in each of the following three configurations. The dimension is $n = 10{,}000$ and there are $k = 100$, $k = 300$ and $k = 500$ nonzero signals. The nonzero signals are selected, by simulation uniformly between $-3$ and $3$, independently in each simulation.

We compare the performance of our procedure with a hard threshold Strong Oracle, whose loss per sample $i$:

$$\min_C \sum_i (Y_i - \mu_i)^2 I(|Y_i| > C) + (0 - \mu_i)^2 I(|Y_i| < C).$$

Thus the Strong Oracle applies the best hard threshold per realization. The entries in Table 2 are based on the average of 50 simulations.

We see that the Strong Oracle, dominates our procedure in the very sparse case where $k = 100$. Our procedure dominates in the less sparse cases.

In the following Table 3, we report on the behavior of our procedure, based on 50 simulations in each of the following configurations. The dimension is $n = 100{,}000$ with $k$ nonzero signals, and each has the value 4. The simulations are performed for $k = 500$, $k = 1000$ and $k = 5000$. The comparison is again with a Strong Oracle. In our procedure we let $v = 1.05$. Note that our procedure dominates the SO in each case. Our procedure thus appears more advantageous as the dimension increases and there are more observations available to estimate $\delta^{*\mu}$.

Table 2
Risk of $\tilde{\delta}_{1.1}$ compared with that of a Strong Oracle; $n = 10{,}000$

|  | $\tilde{\delta}_{1.1}$ | SO |
|---|---|---|
| $k = 100$ | 306 | 295 |
| $k = 300$ | 748 | 866 |
| $k = 500$ | 1134 | 1430 |

TABLE 3
*Risk of $\tilde{\delta}_{1.05}$ compared
with that of a Strong Oracle; $n = 100{,}000$*

|  | $\tilde{\delta}_{1.05}$ | SO |
|---|---|---|
| $k = 500$ | 2410 | 3335 |
| $k = 1000$ | 3810 | 5576 |
| $k = 5000$ | 10,400 | 16,994 |

## APPENDIX

PROOF OF LEMMA 2. Let $r_{v,n}$ be the risk corresponding to $\delta_v^{*\boldsymbol{\mu}}(u)$ when applied on an independent sample $U_i^n \sim N(\mu_i, v)$, and let $r_{1,n}$ be the risk of $\delta_1^{*\boldsymbol{\mu}}(y)$ when applied on an independent sample $Y_i^n \sim N(\mu_i, 1)$, $i = 1, \ldots, n$. Then,

$$E_{\boldsymbol{\mu}^n} \sum (\delta_v^{*\boldsymbol{\mu}}(U_i^n) - \mu_i^n)^2 = r_{v,n}, \tag{37}$$

$$E_{\boldsymbol{\mu}^n} \sum (\delta^{*\boldsymbol{\mu}}(Y_i^n) - \mu_i^n)^2 = r_{1,n}. \tag{38}$$

We will omit the superscript $n$ in the following.

Obviously, $r_{v,n} > r_{1,n}$, since the experiment $Y_i, i = 1, \ldots, n$ dominates the experiment $U_i, i = 1, \ldots, n$ in terms of comparison of experiments. We will first show that for $v = 1 + (1/d_n)$ where $\log(n) = o(d_n)$,

$$r_{1,n}/r_{v,n} \to 1. \tag{39}$$

Let $\phi(u, \mu_i; s^2)$, denote the normal density with variance $s^2$ and mean $\mu_i$. For every $\varepsilon' > 0$, we have

$$r_{v,n} = \sum_i E_{\mu_i}(\delta_v^{*\boldsymbol{\mu}}(U_i) - \mu_i)^2 \leq \sum_i E_{\mu_i}(\delta_1^{*\boldsymbol{\mu}}(U_i) - \mu_i)^2$$

$$= \sum_i \int (\delta_1^{*\boldsymbol{\mu}}(u) - \mu_i)^2 \phi(u, \mu_i; v)\, du$$

$$= \sum_i \int (\delta_1^{*\boldsymbol{\mu}}(u) - \mu_i)^2 \phi(u, \mu_i; 1) \frac{\phi(u, \mu_i; v)}{\phi(u, \mu_i; 1)}\, du \tag{40}$$

$$= (1 + o(1)) \times \sum_i E_{\mu_i}(\delta_1^{*\boldsymbol{\mu}}(Y_i) - \mu_i)^2 + o(n^{\varepsilon'}) \tag{41}$$

$$= (1 + o(1)) \times r_{1,n} + o(n^{\varepsilon'}). \tag{42}$$

Equation (41) is implied as follows. When $d_n/\log(n) \to \infty$, then, for each summand $i$ in (40), the ratio of the densities approaches 1 uniformly on the range where $|u - \mu_i| < \sqrt{K \log(n)}$ for any $K$. It is also easy to check for each summand, that, for large enough $K$, the integral over $|u - \mu_i| >$ $\sqrt{K \log(n)}$ can be made of the order $o(n^{(\varepsilon'-1)})$ for any $\varepsilon'$. This may be seen since $|\delta^{*\boldsymbol{\mu}}|$ and $|\mu_i|$ are bounded by $C_n$, while by choosing $K$ large enough $n \times (2C_n)^2 \times P(|Y_i - \mu_i| > \sqrt{K \log(n)})$ can be made of order $o(n^{\varepsilon'})$.

By (30), letting $\varepsilon' < \varepsilon$, (42) implies $\limsup r_{v,n}/r_{1,n} \leq 1$. This completes the proof of (39), since, as mentioned, $r_{1,n} \leq r_{v,n}$.

Similarly to the above, we write

$$\text{(43)} \quad \sum_i E_{\mu_i}(\delta_v^{*\boldsymbol{\mu}}(Y_i) - \mu_i)^2 = \sum_i \int (\delta_v^{*\boldsymbol{\mu}}(t) - \mu_i)^2 \phi(t, \mu_i; 1)\, dt$$

$$= \sum_i \int (\delta_v^{*\boldsymbol{\mu}}(t) - \mu_i)^2 \frac{\phi(t, \mu_i, 1)}{\phi(t, \mu_i; v)} \phi(t, \mu_i; v)\, dt.$$

An argument similar to the above (yet easier) implies that for $d_n/\log(n) \to \infty$ we have

$$\text{(44)} \quad \frac{r_{v,n}}{E_{\boldsymbol{\mu}} \sum (\delta_v^{*\boldsymbol{\mu}}(Y_i) - \mu_i)^2} \to 1.$$

Lemma 2 now follows from (39) and (44).

Note that Lemma 2 would follow along the same lines if we assume in (30) the weaker condition $R(\boldsymbol{\mu}^n, \Delta^{\boldsymbol{\mu}^n}) = O(1)$ (i.e., under our notion of an extremely sparse setup). □

PROOF OF LEMMA 3. In order to motivate the expression in (46), bellow, we begin by comparing the performance of $\delta_v^{*\boldsymbol{\mu}}$ and $\hat{\delta}_v$ when applied on a set of new independent observations $\tilde{Y}_i^n \sim N(\mu_i^n, 1)$, $i = 1, \ldots, n$, which are also independent of the set $Y_i^n$, $i = 1, \ldots, n$, that was used to obtain the estimate $\hat{\delta}_v$. We will omit the superscript $n$ in the following. Thus, we first show that

$$\text{(45)} \quad E_{\boldsymbol{\mu}} \sum_i (\delta_v^{*\boldsymbol{\mu}}(\tilde{Y}_i) - \hat{\delta}_v(\tilde{Y}_i))^2 = o(n^\varepsilon)$$

for every $\varepsilon > 0$.

Observe that

$$\text{(46)} \quad E_{\boldsymbol{\mu}} \sum_i (\delta_v^{*\boldsymbol{\mu}}(\tilde{Y}_i) - \hat{\delta}_v(\tilde{Y}_i))^2 = E \sum_i \int [\delta_v^{*\boldsymbol{\mu}}(y) - \hat{\delta}_v(y)]^2 \phi(y - \mu_i)\, dy$$

$$= n \int_{-\infty}^{\infty} E[(\delta_v^{*\boldsymbol{\mu}}(y) - \hat{\delta}_v(y))^2] g_{G,1}^*(y)\, dy.$$

Here, $g_{G,1}^*(y) = \frac{1}{n} \sum_i \phi(y - \mu_i)$ denotes the mixture density in the compound decision setup, where $G$ corresponds to the empirical distribution of $(\mu_1, \ldots, \mu_n)$ (see Remark 1, Section 2). In fact, $G \equiv G_n^{\boldsymbol{\mu}}$, but the dependence of $G$ on $n$ and $\boldsymbol{\mu}$ is suppressed in the notation.

The outline of the proof, that (46) is of order $o(n^\varepsilon)$ for every $\varepsilon > 0$, is as follows. Let $0 < \varepsilon' < \varepsilon$. Let $\mathcal{R}$ consist of all point $y$ that satisfy:

(i) $-C'_n < y < C'_n$, where $C'_n = (\log(n) + C_n)$
and,
(ii) $g^*_{G,1}(y) > n^{\varepsilon'-1}$ for some $0 < \varepsilon' < \varepsilon$.

We then show that, uniformly for $y_0 \in \mathcal{R}$,

(47) $$E[(\delta^{*\boldsymbol{\mu}}_v(y_0) - \hat{\delta}_v(y_0))]^2 = o\left(\frac{n^{\varepsilon'}}{ng^*_{G,1}(y_0)}\right).$$

Once (47) has been verified, the proof of (45) can be completed, since

$$n \int_{-\infty}^{\infty} E[(\delta^{*\boldsymbol{\mu}}_v(y) - \hat{\delta}_v(y))^2] g^*_{G,1}(y) \, dy$$

$$= n \int_{\mathcal{R}^c \cap [-C'_n, C'_n]} E[(\delta^{*\boldsymbol{\mu}}_v(y) - \hat{\delta}_v(y))]^2 g^*_{G,1}(y) \, dy$$

(48)
$$+ n \int_{y \notin [-C'_n, C'_n]} [(\delta^{*\boldsymbol{\mu}}_v(y) - \hat{\delta}_v(y))]^2 g^*_{G,1}(y) \, dy$$

$$+ n \int_{\mathcal{R}} E[(\delta^{*\boldsymbol{\mu}}_v(y) - \hat{\delta}_v(y))]^2 g^*_{G,1}(y) \, dy$$

$$= o(n^\varepsilon) + \left(n \int_{\mathcal{R}} \frac{o(n^{\varepsilon'})}{ng^*_{G,1}(y)} g^*_{G,1}(y) \, dy\right)$$

(49)
$$= o(n^\varepsilon) + o(C'_n n^{\varepsilon'}) = o(n^\varepsilon).$$

In the above, we use the exponential tail of the normal distribution, the truncation, and the upper bound $C_n$ for the elements of $\boldsymbol{\mu}^n$.

We elaborate now on the derivation of (47). The (nontruncated) version of our estimator equals

(50) $$\hat{\delta}_v(y) = y + \frac{\hat{g}^{*\prime}_h(y)}{\hat{g}^*_h(y)},$$

where $\hat{g}^*_h$ and $\hat{g}^{*\prime}_h$ are kernel density estimators, based on $Y_1, \ldots, Y_n$, of $g^*_{G,1}$ and $g^{*\prime}_{G,1}$ with bandwidth $h = \sqrt{v-1} \equiv \sqrt{1/d_n}$.

Recall that

(51) $$\delta^{*\boldsymbol{\mu}}_v(y) = y + v \frac{g^{*\prime}_{G,v}(y)}{g^*_{G,v}(y)}.$$

We now write the right-hand side of (50) as

(52) $$\hat{\delta}_v(y) = y + \frac{g^{*\prime}_{G,v}(y) + R_1}{g^*_{G,v}(y) + R_2}.$$

Here, the random variables $R_1$ and $R_2$ are implicitly defined, by comparing numerators and denominators of (52) and (50), respectively.

Note that, by Lemma 1,

(53) $E(R_i) = 0, \quad i = 1, 2,$

(54)
$$\begin{aligned}
E[\delta_v^{*\boldsymbol{\mu}}(y_0) - \hat{\delta}_v(y_0)]^2 \\
&= O\bigg(E\bigg(\frac{R_1}{g_{G,v}^*(y_0) + R_2}\bigg)^2 + E\bigg(\frac{g_{G,v}^{*\prime}(y_0)R_2}{(g_{G,v}^*(y_0))^2 + g_{G,v}^*(y_0)R_2}\bigg)^2\bigg) \\
&= O\bigg(E\bigg(\frac{R_1}{g_{G,v}^*(y_0) + R_2}\bigg)^2 + E\bigg(\frac{C_n R_2}{g_{G,v}^*(y_0) + R_2}\bigg)^2\bigg) \\
&= O\bigg(E\bigg(\frac{R_1}{g_{G,1}^*(y_0) + R_2}\bigg)^2 + E\bigg(\frac{C_n R_2}{g_{G,1}^*(y_0) + R_2}\bigg)^2\bigg).
\end{aligned}$$

For the last equality, we use the fact that $g_{G,1}^*(y)/g_{G,v}^*(y)$ is bounded when $v > 1$; for the previous one, we use the fact that $g_{G,v}^{*\prime}(y_0)/g_{G,v}^*(y_0) = O(C_n)$, uniformly for $y_0 \in \mathcal{R}$.

The assertion $E[(\delta_v^{*\boldsymbol{\mu}}(y_0) - \hat{\delta}_v(y_0))]^2 = \frac{o(n^{\varepsilon'})}{ng_{G,1}^*(y_0)}$, will be implied by computing the variances of $R_i, i = 1, 2$, and by controlling the moderate deviation of $R_2$, as in what follows.

The variances of $R_1$ and $R_2$ equal to the variances of the corresponding kernel density estimators in (50), of the density and its derivative. It may be checked, from (18) and (19), that

(55) $$\operatorname{var}(R_i) = O\bigg(\frac{(C_n d_n)^2 g_{G,1}^*(y_0)}{n}\bigg) = \frac{o(n^{\tilde{\varepsilon}}) g_{G,1}^*(y_0)}{n}$$

for every $0 < \tilde{\varepsilon} < \varepsilon'$.

Applying Bernstein's inequality [see, e.g., van der Vaart and Wellner (1996), page 103] we obtain

(56) $$P(R_2 < -0.5 g_{G,1}^*(y_0)) < 1/4C_n^2.$$

Since $[\hat{\delta}_v(y_0) - \delta^{*\boldsymbol{\mu}}(y_0)]^2 < 4C_n^2$ by truncation, (47) follows when incorporating the above computed values of the second moments of $R_i$, $i = 1, 2$ into the numerator of (54) and controlling its denominator by (56).

It remains to show how to modify the proof of (45) in order to conclude $E_{\boldsymbol{\mu}} \sum (\delta_v^{*\boldsymbol{\mu}}(Y_i) - \hat{\delta}_v(Y_i))^2 = o(n^\varepsilon)$ for every $\varepsilon > 0$. We briefly explain it in the following.

Let $\hat{\delta}_v^{(i)}$ be our estimator for $\delta^{*\boldsymbol{\mu}}$ based on $Y_j$, $j = 1, \ldots, n, j \neq i$. We now write

(57)
$$\begin{aligned}
E_{\boldsymbol{\mu}} \sum (\delta_v^{*\boldsymbol{\mu}}(Y_i) - \hat{\delta}_v(Y_i))^2 \\
= E_{\boldsymbol{\mu}} \sum_i (\delta_v^{*\boldsymbol{\mu}}(Y_i) - \hat{\delta}_v^{(i)}(Y_i) + \hat{\delta}_v^{(i)}(Y_i) - \hat{\delta}_v(Y_i))^2.
\end{aligned}$$

We now show that
$$E_{\boldsymbol{\mu}} \sum_i (\hat{\delta}_v(Y_i) - \hat{\delta}_v^{(i)}(Y_i))^2 = o(n^\varepsilon) \tag{58}$$

for every $\varepsilon > 0$. This follows by arguments similar to the ones presented in the first part of our lemma. Specifically, first note that

$$\begin{aligned}
E_{\boldsymbol{\mu}} &\sum_i (\hat{\delta}_v(Y_i) - \hat{\delta}_v^{(i)}(Y_i))^2 \\
&= E_{\boldsymbol{\mu}} \sum_i (\hat{\delta}_v(Y_i) - \hat{\delta}_v^{(i)}(Y_i))^2 \times I(Y_i \in \mathcal{R}) + o(n^\varepsilon) \\
&= o(n^\varepsilon) + E_{\boldsymbol{\mu}} \sum \frac{o(n^{\varepsilon'})}{n\hat{g}^*(Y_i)}.
\end{aligned} \tag{59}$$

Now, taking a dense enough grid in the region $\mathcal{R}$ [note the derivative of $g^*_{G,1}$ in that region is bounded by $O(C_n)$], and applying Bernstein's inequality coupled with Bonferroni, yields

$$P_{\boldsymbol{\mu}}\left(\sup_{y \in \mathcal{R}} \frac{\hat{g}^*(y)}{g^*(y)} < \frac{1}{2}\right) = o\left(\frac{1}{nC_n^2}\right). \tag{60}$$

By (59) and (60) and the truncation, we obtain

$$E_{\boldsymbol{\mu}} \sum \frac{o(n^{\varepsilon'})}{n\hat{g}^*(Y_i)} = o(n^\varepsilon) + \int_{\mathcal{R}} \frac{2 \times o(n^{\varepsilon'})}{g^*(y)} g^*(y)\, dy = o(n^\varepsilon). \tag{61}$$

The above involves interchanging summation and integration.

We then note that
$$\begin{aligned}
E_{\boldsymbol{\mu}} &\sum_i (\delta_v^{*\boldsymbol{\mu}}(Y_i) - \hat{\delta}_v^{(i)}(Y_i))^2 \\
&= E_{\boldsymbol{\mu}} \sum_i \int (\delta_v^{*\boldsymbol{\mu}}(y) - \hat{\delta}_v^{(i)}(y))^2 \phi(y - \mu_i)\, dy \\
&= E_{\boldsymbol{\mu}} \sum_i \int (\delta_v^{*\boldsymbol{\mu}}(y) - \hat{\delta}_v(y) + \hat{\delta}_v(y) - \hat{\delta}_v^{(i)}(y))^2 \phi(y - \mu_i)\, dy \\
&= o(n^{\varepsilon'})
\end{aligned} \tag{62}$$

for every $\varepsilon > 0$.

Obtaining the last equality involves evaluating
$$E_{\boldsymbol{\mu}} \sum_i \int E(\delta_v^{*\boldsymbol{\mu}}(y) - \hat{\delta}_v(y))^2 \phi(y - \mu_i)\, dy$$

as in (46), and
$$\int \sum_i E(\hat{\delta}_v(y) - \hat{\delta}_v^{(i)}(y))^2 \phi(y - \mu_i)\, dy$$

similarly to (58). This completes the proof. □

**Acknowledgments.** Helpful discussions with Jayanta K. Ghosh, Junyong Park, Ya'acov Ritov and Cun-Hui Zhang, are gratefully acknowledged. We are also grateful to the referee and Associate Editor.

DEPARTMENT OF STATISTICS
THE WHARTON SCHOOL
UNIVERSITY OF PENNSYLVANIA
400 JON M. HUNTSMAN HALL
3730 LOCUST WALK, PHILADELPHIA
PENNSYLVANIA 19104
USA
E-MAIL: lbrown@wharton.upenn.edu

DEPARTMENT OF STATISTICS
DUKE UNIVERSITY
BOX 90251, DURHAM
NORTH CAROLINA, 27708-0251
USA
E-MAIL: eitan.greenshtein@gmail.com